\documentclass[12pt]{article}

\usepackage{amsmath,amsthm,amssymb, amscd, amsfonts}

\begin{document}
\newcommand{\nc}{\newcommand}
\nc{\rnc}{\renewcommand} \nc{\nt}{\newtheorem}


\nc{\TitleAuthor}[2]{\nc{\Tt}{#1}%
    \nc{\At}{#2}%
    \maketitle%
}

\nt{theo}{Theorem}[section] \nt{coro}[theo]{Corollary}
\nt{lemm}[theo]{Lemma} \nt{prop}[theo]{Proposition}
\nt{ques}[theo]{Question} \nt{defi}[theo]{Definition}

\theoremstyle{remark}

\nt{axio}{Axiom}
\newtheorem{conj}{Conjecture}
\renewcommand{\theconj}{}

\nc{\Rm}[2]{\newenvironment{#1}{\begin{trivlist}%
    \addtocounter{theo}{1}%
    \item[]#2~{\bf \thetheo.}}{\end{trivlist}}
}

\Rm{rema}{\emph{Remark}} \Rm{nota}{\emph{Notation}} \Rm{exam}{{\bf
Example}}

\newenvironment{pf}{\begin{trivlist}\item[]{\bf Proof.}~}%
    {~\hfill~\fbox{}\end{trivlist}
}

\nc{\thlabel}[1]{\label{theo:#1}}
\nc{\thref}[1]{Theorem~\ref{theo:#1}}
\nc{\selabel}[1]{\label{sect:#1}}
\nc{\seref}[1]{Section~\ref{sect:#1}}
\nc{\lelabel}[1]{\label{lemm:#1}}
\nc{\leref}[1]{Lemma~\ref{lemm:#1}}
\nc{\prlabel}[1]{\label{prop:#1}}
\nc{\prref}[1]{Proposition~\ref{prop:#1}}
\nc{\colabel}[1]{\label{coro:#1}}
\nc{\coref}[1]{Corollary~\ref{coro:#1}}
\nc{\exlabel}[1]{\label{exam:#1}}
\nc{\exref}[1]{Example~\ref{exam:#1}}
\nc{\delabel}[1]{\label{defi:#1}}
\nc{\deref}[1]{Definition~\ref{defi:#1}}
\nc{\eqlabel}[1]{\label{equa:#1}}

\nc{\Hom}{\operatorname{Hom}} \nc{\Mor}{\operatorname{Mor}}
\nc{\Aut}{\operatorname{Aut}} \nc{\Ann}{\operatorname{Ann}}
\nc{\Ker}{\operatorname{Ker}} \nc{\Trace}{\operatorname{Trace}}
\nc{\Char}{\operatorname{Char}} \nc{\Mod}{\operatorname{Mod}}
\nc{\End}{\operatorname{End}} \nc{\Spec}{\operatorname{Spec}}
\nc{\Span}{\operatorname{Span}} \nc{\sgn}{\operatorname{sgn}}
\nc{\Id}{\operatorname{Id}} \nc{\Com}{\operatorname{Com}}

\nc{\nd}{\mbox{$\not|$}} 
\nc{\nci}{\mbox{$\not\subseteq$}}
\nc{\scontainin}{\mbox{$\mbox{}\subseteq\hspace{-1.5ex}\raisebox{-.5ex}{$_\prime$}\hspace*{1.5ex}$}}


\nc{\R}{{\sf R\hspace*{-0.9ex}\rule{0.15ex}%
    {1.5ex}\hspace*{0.9ex}}}
\nc{\N}{{\sf N\hspace*{-1.0ex}\rule{0.15ex}%
    {1.3ex}\hspace*{1.0ex}}}
\nc{\Q}{{\sf Q\hspace*{-1.1ex}\rule{0.15ex}%
       {1.5ex}\hspace*{1.1ex}}}
\nc{\C}{{\sf C\hspace*{-0.9ex}\rule{0.15ex}%
    {1.3ex}\hspace*{0.9ex}}}
\nc{\Z}{\mbox{${\sf Z}\!\!{\sf Z}$}}


\newcommand{\gd}{\delta}
\newcommand{\sub}{\subset}
\newcommand{\cntd}{\subseteq}
\newcommand{\go}{\omega}
\newcommand{\Pa}{P_{a^\nu,1}(U)}
\newcommand{\fx}{f(x)}
\newcommand{\fy}{f(y)}
\newcommand{\gD}{\Delta}
\newcommand{\gl}{\lambda}
\newcommand{\half}{\frac{1}{2}}
\newcommand{\ga}{\alpha}
\newcommand{\gb}{\beta}
\newcommand{\gga}{\gamma}
\newcommand{\ul}{\underline}
\newcommand{\ol}{\overline}
\newcommand{\Lrraro}{\Longrightarrow}
\newcommand{\equi}{\Longleftrightarrow}
\newcommand{\gt}{\theta}
\newcommand{\op}{\oplus}
\newcommand{\Op}{\bigoplus}
\newcommand{\CR}{{\cal R}}
\newcommand{\tr}{\bigtriangleup}
\newcommand{\grr}{\omega_1}
\newcommand{\ben}{\begin{enumerate}}
\newcommand{\een}{\end{enumerate}}
\newcommand{\ndiv}{\not\mid}
\newcommand{\bab}{\bowtie}
\newcommand{\hal}{\leftharpoonup}
\newcommand{\har}{\rightharpoonup}
\newcommand{\ot}{\otimes}
\newcommand{\OT}{\bigotimes}
\newcommand{\bwe}{\bigwedge}
\newcommand{\gep}{\varepsilon}
\newcommand{\gs}{\sigma}
\newcommand{\OO}{_{(1)}}
\newcommand{\TT}{_{(2)}}
\newcommand{\FF}{_{(3)}}
\newcommand{\minus}{^{-1}}
\newcommand{\CV}{\cal V}
\newcommand{\CVs}{\cal{V}_s}
\newcommand{\slp}{U_q(sl_2)'}
\newcommand{\olp}{O_q(SL_2)'}
\newcommand{\slq}{U_q(sl_n)}
\newcommand{\olq}{O_q(SL_n)}
\newcommand{\un}{U_q(sl_n)'}
\newcommand{\on}{O_q(SL_n)'}
\newcommand{\ct}{\centerline}
\newcommand{\bs}{\bigskip}
\newcommand{\qua}{\rm quasitriangular}
\newcommand{\ms}{\medskip}
\newcommand{\noin}{\noindent}
\newcommand{\raro}{\rightarrow}
\newcommand{\alg}{{\rm Alg}}
\newcommand{\rcom}{{\cal M}^H}
\newcommand{\lcom}{\,^H{\cal M}}
\newcommand{\rmod}{\,_R{\cal M}}
\newcommand{\qtilde}{{\tilde Q^n_{\gs}}}
\nc{\e}{\overline{E}} \nc{\K}{\overline{K}}
 \nc{\lyd}{\,_H^H{\cal YD}} \nc{\mch}{\,_H^H{\cal
M}} \nc{\hb}{\tilde{H}} \nc{\frt}{f_{R^{\tau}}}
\def\newtheorems{\newtheorem{theorem}{Theorem}[subsection]
         \newtheorem{cor}[theorem]{Corollary}
         \newtheorem{prop}[theorem]{Proposition}
         \newtheorem{lemma}[theorem]{Lemma}
         \newtheorem{defn}[theorem]{Definition}
         \newtheorem{Theorem}{Theorem}[section]
         \newtheorem{Corollary}[Theorem]{Corollary}
         \newtheorem{Proposition}[Theorem]{Proposition}
         \newtheorem{Lemma}[Theorem]{Lemma}
         \newtheorem{Defn}[Theorem]{Definition}
         \newtheorem{Example}[Theorem]{Example}
         \newtheorem{Remark}[Theorem]{Remark}
         \newtheorem{claim}[theorem]{Claim}
         \newtheorem{sublemma}[theorem]{Sublemma}
         \newtheorem{example}[theorem]{Example}
         \newtheorem{remark}[theorem]{Remark}
         \newtheorem{question}[theorem]{Question}
         \newtheorem{conjecture}{Conjecture}[subsection]}

\title{Finite-dimensional pointed Hopf algebras of type $A_n$ related to
the Faddeev-Reshetikhin-Takhtajan $U(R)$ construction}
\author{Jacob Towber$^1$
\\ Department of Mathematics, Statistics and Computer Science\\
University of Illinois at Chicago, Chicago, Illinois\\
\and Sara Westreich$^1$
\\Interdisciplinary Department of the Social
Sciences
\\Bar-Ilan University,  Ramat-Gan, Israel}
\footnotetext[1]{The authors wish to thank UIC, where some of the
work was done, for hospitality.}

 \maketitle

\begin{abstract}
Two "quantum enveloping algebras", here denoted by $U(R)$ and
$U^{\sim}(R)$, are associated in \cite{frt88} and \cite{frt89} to
any Yang-Baxter operator R. The latter is only a bialgebra, in
general; the former is a Hopf algebra.

In this paper, we study the pointed Hopf algebras $U(R_Q)$, where
$R_Q$ is the Yang-Baxter operator associated with the
multi-parameter deformation of $GL_n$ supplied in \cite{ast};  cf
also \cite{s,re}.

Some earlier results concerning these Hopf algebras $U(R_Q)$ were
obtained in \cite{t,clmt,cm}; a related (but different) Hopf
algebra was studied in \cite{dp}.

The main new results obtained here concerning these quantum
enveloping algebras are: 1)We list, in an extremely explicit form,
those quantum enveloping algebras $U(R_Q)$ which are
finite-dimensional---let ${\mathcal U}$ denote the collection of
these. 2)We verify that the pointed Hopf algebras in ${\mathcal
U}$ are quasitriangular and of Cartan type $A_n$ in the sense of
Andruskiewich-Schneider.
 3)We show that every $U(R_Q)$ is a Hopf quotient of a double cross-product
(hence, as asserted in 2), is quasitriangular if
finite-dimensional.) 4) CAUTION: These Hopf algebras are NOT
always
 cocycle twists of the standard 1-parameter deformation. This somewhat surprising
 fact is an immediate
 consequence of the data furnished here---
clearly a cocycle twist will not convert an infinite-dimensional
Hopf algebra to a finite-dimensional one! Furthermore, these Hopf
algebras in ${\mathcal U}$ are (it is proved) not all cocycle
twists of each other. 5)We discuss also the case when the quantum
determinant is central in $A(R_Q)$, so it makes sense to speak of
a $Q$-deformation of the special linear group.
\end{abstract}

\section*{Inroduction}
Throughout the remainder of this paper, $k$ will denote an
algebraically closed ground-field of characteristic $0$.

Let $G$ be an affine algebraic group over $k$; we then associate
to $G$ in the usual way two $k$-Hopf algebras : $A(G)$, whose
elements are representative functions on $G$, and $U(G)$, whose
underlying $k$-algebra is the enveloping algebra of the Lie
algebra of $G$. Those Hopf algebras certified by workers in the
field as being ``quantum groups'', fall into two main classes:
those `deforming' the type $U(G)$, the ``quantum enveloping
algebras'', and those `deforming' the type $A(G)$, which will here
be called ``quantum groups''.

Perhaps the earliest systematic construction of infinite families
of these two types of Hopf algebras, was furnished by the seminal
work of Faddeev, Reshetikhin and Takhtejan (\cite{frt88,frt89}).
The starting point of their marvellous construction may be taken
to be a $K$-linear transformation
$$R:V\otimes V \mapsto V\otimes V $$
where $V$ is a finite-dimensional vector-space over $K$, and $R$
satisfies the Yang-Baxter condition, say in the form
\begin{equation}\label{yb}(R\otimes I_V)\circ (I_V\otimes R)\circ (R\otimes I_V)=
(I_V\otimes R)\circ (R\otimes I_V)(I_V\otimes R) \end{equation}

Given such an $R$, there is associated in \cite{frt88} a Hopf
algebra $A(R)$, first proved in \cite{lt} to have the property
(which seems characteristic for ``quantum groups'') that its
finite-dimensional comodules form in a natural way a braided
monoidal category (as defined in \cite{js}). Moreover, the paper
\cite{frt88} goes on to construct inside $(A(R))^\circ$ a
bialgebra, which here will be denoted by $U^{\sim}(R)$; this is in
general only a bialgebra, but not a Hopf algebra. This result was
improved by Faddeev, Reshetikhin and Takhtejan in a later
paper\cite{frt89}, where they construct inside $(A(R))^\circ$ a
bialgebra---which will here be denoted by $U(R)$---properly
containing the earlier construction $U^{\sim}(R)$, and where they
show that this larger bialgebra $U(R)$  indeed has in a natural
way the structure of a Hopf algebra.

The most usual applications of quantum groups, have involved the
$1$-parameter deformations of quantum enveloping algebras, given
by the Drin'feld-Jimbo-Lusztig construction. Thus, there was some
interest aroused in the early '90s, by the construction of
$[{N\choose 2}+1]$-parameter deformations of $A(GL_n)$ and
$U(gl_n),$  as given in \cite{ast,re,s}. This work involved the
construction of a solution $R=R_Q$ to \eqref{yb}, where the
multiparameter
$$Q=\{r,q_{i,j}:1\leq i<j\leq N\} $$
is made up of $[{N\choose 2}+1]$ non-zero elements $r,q_{i,j}$ in
$k$.

 We are thus led to study the quantum enveloping algebra $U(R_Q)$
(furnished by applying to $R_Q$ the constructions of
Faddeev,Reshitikhin and Takhtajan discussed above.) This was first
done in complete generality  in the paper \cite{t}, while in the
later-appearing papers \cite{clmt} and \cite{cm}, the proofs in
\cite{t} were substantially simplified, at the expense of
requiring the parameter $r$ in $Q$, not to be a root of $1$---this
additional assumption is however {\it never} satisfied for the
cases to be considered in the present paper, as will be explained
below.

\medskip
\medskip
Actually, the Hopf algebra studied in \cite{t}, only coincides
with the \cite{frt89} construction $U(Q)$ when $r \neq 1$;when
$r=1$ it represents a `closure' as $r\rightarrow 1$ whose further
study we reserve for a later paper.
\medskip
For the sake of completeness, let us mention yet another
incarnation of $U(R_Q)$, constructed (in work appearing prior to
\cite{t} ) by Dobrev and Parashar, in \cite{dp}. Their Hopf
algebra (whose definition utilizes $R_Q$ in an interesting way
which is outside the scope of the present paper) is never
finite-dimensional, hence is not directly relevant to our present
purposes. In the addendum to \cite{dp} they give a multiparameter
$Q$-deformation of $U(sl_n)$---which also is never
finite-dimensional.

\medskip There has recently been an interest in the study of pointed finite-dimensional
 Hopf algebras, stimulated by the remarkable results obtained in this
direction by Andruskiewich and Schneider(\cite{as,as2}). The
purpose of the present paper is to study in this light, the rather
ancient work on $U(R_Q)$ discussed above. Since we are thus here
interested only in finite-dimensional Hopf algebras, only the
construction presented in \cite{t} will be relevant.
\medskip
It will be proved below that, for $r\neq 1$,  $U(R_Q)$ is
finite-dimensional, if and only if:
\smallskip
 \noindent (*)\ \ $r$ and each $q_{i,j}$ is a root of unity.

 (From this it is easily deduced that, as asserted above, none of the versions
of $U(R_Q)$ constructed in \cite{re,s,dp,clmt,cm}is ever
finite-dimensional.) Also, it will be proved below that, when (*)
holds, the pointed Hopf algebra $U(R_Q)$ is quasi-triangular, and
has Andruskiewich-Schneider Cartan matrix of type (as it `should'
be!)$A_n$.
 These Hopf algebras are not all twistings of
each other though. Indeed, we exhibit in Theorem 2.2 families of
finite-dimensional Hopf algebras arising from $U(R_Q)$ with
non-isomorphic groups of group-like elements, and thus these Hopf
algebras can not be obtained from each other by twists. Although
not all pointed Hopf algebras of type $A_n$ can be realized as
$U(R_Q),$ the particular new family described in
\cite[Ex.7.27]{as} can be so realized.

The structure of the group of group-like elements of $U(R_Q)$
strongly depends on the determinant element of $A(R_Q).$ Recall
that the bialgebra $O_q(M_n(\mathcal{C}))$ can be obtained as
$A(R_Q)$ for a special choice of the parameters in $Q$ and that
its determinant group-like element is central. The Hopf algebra
$U_q(gl_n)$ and  Lusztig's finite-dimensional Hopf algebra
$u_q(gl_n)$ are both obtained as the corresponding $U(R_Q),$
depending on whether or not $q$ is a root of unity, while
$U_q(sl_n)$ and $u_q(sl_n)$ are Hopf subalgebras respectively. We
show in \thref{family} how this situation is generalized for a
Hopf algebra of the form $U(R_Q)$ where $A(R_Q)$ admits a central
determinant. Example 2.5 is then a special examples of a (new)
Hopf algebra of this type. We construct a Hopf algebras of type
$A_3$ with a group of group-like elements generated by one element
(while the group of group-like elements of $u^{\ge 0}_q(sl_3)$ is
generated by two elements).

In the second part we prove in \thref{double} that  $U(R_Q)$ is
always a Hopf quotient of a double crossproduct of its "$\ge 0$"
and "$\le 0$" parts. This implies in particular that $U(R_Q)$ is
quasitriangular when it is finite dimensional.

\section{Preliminaries}
Throughout we assume that the base field $k$ is algebraically
closed of characteristic $0.$\medskip

\noin\ul{gr$H$ as a biproduct:}\medskip

In what follows we give is a brief overview of this subject based
on one of the many possible references (see for example
\cite{as,as2}).

Let $H$ be a pointed  Hopf algebra over an algebraically closed
field of characteristic $0$, let $H_n,\,n\geq 0$ denote the
coradical filtration of $H$ and set $H_{-1}=k.$ Let $G=G(H)=H_0$
 and let
$${\rm gr}\,H=\underset{n\geq 0}\oplus{\rm gr}\,H(n)$$ where
${\rm gr}\,H(n)=H_n/H_{n-1}$ for all $n\geq 0.$ Then ${\rm gr}\,H$
is a graded Hopf algebra. There is a Hopf algebra projection
$\pi:{\rm gr}\,H\rightarrow {\rm gr}\,H(0)=kG$ and a Hopf algebra
injection $i:kG\rightarrow {\rm gr}\,H.$ By \cite{ra} this implies
that we have a biproduct
$${\rm gr}\,H\cong R\#kG$$
where $R=\{x\in {\rm gr}\,H\,|\,(id\ot \pi)\circ\gD(x)=x\ot 1\}$
is the algebra of the coinvariants of the induced $H$-coaction.

It is known that $R$ is a graded braided Hopf algebra in the
category of left Yetter-Drinfeld modules over $kG.$ The action of
$kG$ on $R$ is given by the adjoint action of the group and the
coaction is given by $\rho=(\pi\ot id)\circ\gD.$ The original Hopf
algebra is then a {\bf lifting} of $R\#kG.$

By \cite{ra} there is a coalgebra projection $\Pi:R\#kG\rightarrow
R$ given by
$$\Pi=id*(i\circ S\circ \pi)$$
The following lemma follows directly from the definition of $\Pi.$
We include it here for completeness.
\begin{lemm}\lelabel{primitive}
All group-like elements of $R\#kG$ are mapped by $\Pi$ to $1.$ A
skew primitive element $x$ such that $\gD(x)=g'\ot x+x\ot
g,\;g,g'\in G$  is mapped to $xg\minus$ which is a primitive
element of $R.$
\end{lemm}

The vector space $V=P(R)$ of primitive elements of $R$ is a
Yetter-Drinfeld submodule of $R$ with a braiding (called the {\bf
infinitesimal braiding})
$$c:V\ot V\rightarrow V\ot V\quad\text{given by}\quad c(v\ot w)=\sum(v_{-1}\cdot w)\ot v$$
where $\rho(v)=\sum v_{-1}\ot v_0\in H\ot V.$

The {\bf Nichols algebra} of $V,\; B(V),$ is in this case the
subalgebra of $R$ generated by $V.$

If the group $G$ is abelian and $V$ is finite-dimensional then the
braiding is given by a family of scalars $l_{ij}\in k,1\leq
i,j\leq n$ so that
$$c(x_i\ot x_j)=l_{ij}(x_j\ot x_i)$$
where $\{x_1,\dots,x_n\}$ is a basis of $V.$ We say that the
braiding is of {\bf Cartan-FL-type} if there exist $q\ne 1$ so
that for $1\leq i,j\leq n$
$$l_{ij}l_{ji}=q^{d_ia_{ij}}$$
where $(a_{ij})$ is a generalized symmetrizable Cartan matrix with
positive integers $\{d_1,\dots,d_n\} $ so that
$d_ia_{ij}=d_ja_{ji}.$

The Cartan matrix is  invariant under twisting (in the sense of
\cite[\S 2]{as2} which is a variation of Reshetikhin\cite{re}).
More precisely, a twist for a Hopf algebra $H$ is an invertible
element $F\in H\ot H$ which satisfies

$$(\Delta\ot Id)(F)(F\ot 1)=(Id\ot \Delta)(F)(1\ot F)$$
and $$ (\varepsilon\ot Id)(F)=(Id\ot \varepsilon)(F)=1.$$

Given a twist $F$ for $H,$ one can define a new Hopf algebra $H^F$
where $H^F=H$ as an algebras and the coproduct is determined by
$$\Delta^ F(a)=F^{-1}\Delta(a)F,\;S^F(a)=Q^{-1}S(a)Q$$
for every $a\in H,$ where $Q:=m\circ(S\ot Id)(F).$

In this context we are interested in the particular case when
$F\in kG\ot kG$ and the group $G=G(H)$ is commutative.

Let $\widehat{G}$ denote the group of characters of the abelian
group $G$ and let $\sigma$ be a (convolution)-invertible
$2$-cocycle on $\widehat{G}.$  Then $\sigma$ gives rise to a twist
$F\in kG\ot kG.$ It is proved that the infinitesimal braiding of
$H^F$ is of the same Cartan matrix as that of $H.$
\bigskip

\noin\ul{The Hopf algebras $A(R_Q)$ and $U(R_Q)$:}\medskip

The reader is referred to \cite{t} for the full details; we follow
the notations there.

Let $Q=\{r\ne 1,\,p_{i,j}\}_{1\le i<j\le n}$ be  ${n\choose 2}+1$
non-zero elements of $k$ and
$$q_{i,j}:=r/p_{i,j}.$$ Set
\begin{equation}\label{kappa}\kappa_j^i=
\begin{cases}
p_{ij} \quad &i<j,\\
r &i=j,\\
q_{ji}=r/p_{ji} & i>j
\end{cases} \end{equation}
Let $V$ be a vector space with a basis $\{v_1,\dots,v_n\}$ and let
$R_Q:V\ot V\rightarrow V\ot V$  be the following Yang-Baxter
operator:
$$R_Q(v_i\ot v_j)=\sum_{k,l=1}^nR_{ij}^{kl} v_k\ot v_l$$
where \begin{equation}\label{rijkl} R_{ij}^{kl}=
\begin{cases}
\kappa_j^i & i=l,\,j=k\\
r-1 &i=k,\,j=l,\,i>j\\
0& otherwise
\end{cases}\end{equation}
Let $A(R_Q)$ be the associated bialgebra constructed by the
FRT-construction and described in \cite{ast}. Recall that $A(R_Q)$
is generated as an algebra by $\{T_i^j\}_{1\le i,j\le n}$ so that
\begin{equation}\label{deltaoq}\gD(T_i^j)=\sum_k T_i^k\ot
T_k^j.\end{equation}

 Let $U_Q=U(R_Q)\subset (A(R_Q))^0$ be the
FRT-construction of the $U$-Hopf algebra described in \cite{t}.

We recall that for $1\le i\le n$ the group-like elements
$K_i,L_i\in U_Q$ are the algebra automorphisms defined on the
generators $T_j^l$ of $A(R_Q)$ by
\begin{equation}\label{ki}
 K_i(T_j^l)=\gd_{jl}\kappa_j^i\qquad L_i(T_j^l)=\gd_{jl}(\kappa _i^j)\minus \end{equation}
and $G(U_Q)$ is an abelian group generated by $\{K_i^{\pm
1},L_i^{\pm 1}\}.$

There are skew-primitive elements $E_{i+1}^i\in
U_Q,\;i=1,\dots,n-1$ defined on the generators $T_j^l$ of $A(R_Q)$
by
$$E^i_{i+1}(T_k^l)=\gd_{k,i}\gd_{l,i+1}$$
so that
\begin{equation}\label{deltae}\gD(E_{i+1}^i)=K_{i+1}\ot E_{i+1}^i+E_{i+1}^i\ot K_i.
\end{equation}
There are also  skew-primitive elements $F^{i+1}_i\in U_Q$ defined
by
$$F^{i+1}_i(T_k^l)=\gd_{k,i+1}\gd_{l,i}$$
so that
\begin{equation}\label{deltaf} \gD(F_i^{i+1})=L_i\ot
F_i^{i+1}+F_i^{i+1}\ot L_{i+1}
\end{equation}
 The commuting relations for the $K_i$'s and the
$E_{j+1}^j$'s are given by:
\begin{equation}\label{commute}
\kappa_{j+1}^i K_i E_{j+1}^j=\kappa_j^iE_{j+1}^jK_i
\end{equation}
 The algebra $U_Q$ is generated by
 $\{K_t,\,L_t,\,E_{i+1}^i,\,F_{i}^{i+1}\}_{1\le t\le n,\,1\le i\le n-1}.$

It was proved \cite[Th. 8.5]{t} that $U_Q$ has a PBW basis given
by products of powers of elements $\{E_j^i,F_k^l,g\}_{1\le i< j\le
n,\,1\le k< l\le n}$ where $g\in G(U_Q).$ The elements
$\{E_j^i\}_{i<j}$ can be defined successively starting from the
elements $\{K_t,\,E_{i+1}^i\}$ via the identities in
\cite[(5.14)$^+$-(5.18)$^+$]{t}. Similarly, the elements
$\{F_k^l\}_{k<l}$ can be defined successively starting from the
elements $\{L_t,\,F^{i+1}_i\}$ via the identities in
\cite[(5.14)$^-$-(5.18)$^-$]{t}. If $r$ is a root of unity then
the additional identities
$$(E_i^j)^{e(Q)}=(F_i^j)^{e(Q)}=0$$ hold, where $e(Q)$
is defined in \cite[(8.1)]{t}. Any other relation among the
$E_i^j$'s or among the $F_k^l$'s can be derived from these sets of
identities.

\section{Hopf algebras of type $A_n$ arising from $U_Q$}

In this section we discuss gr$\,U_Q$ and show that it provides new
examples of Hopf algebras  of type $A_n.$

Consider the (Hopf)-subalgebra $B^l\subset U_Q$ generated by :
$$ B^l=<K_j^{\pm 1},E_{i+1}^i,\;j=1,\dots n,\,i=1,\dots,n-1>$$
 Let
$\rm{gr}\,B^l=R\#kG(B^l).$ By \leref{primitive} we have that
\begin{equation}\label{v}
V=Sp_k\{x_i:= E_{i+1}^iK_i\minus,\;i=1,\dots n-1\} \end{equation}
is a subspace of primitive elements contained in $R.$ Observe that
for each $i,$
$$\rho(x_i)=K_i\minus K_{i+1}\ot x_i.$$

Let $SLG$ be the group generated by the $n-1$ group-like elements
\begin{equation}\label{group}
SLG=<\overline{K_i}=K_i\minus K_{i+1},\;i=1,\dots,n-1>
\end{equation}
Then $V$ is a Yetter-Drinfeld module over $kSLG.$ We consider the
Hopf algebra $B(V)\#kSLG\subset {\rm gr}\,B^l.$

Since the Yang-Baxter operator $R_Q$ is related to $GL_n$ one
would expect The following proposition.
\begin{prop}\prlabel{an}
Let $V$ and $SLG$ be as above, then $B(V)\#kSLG$ is of type $A_n.$
\end{prop}
\begin{pf}
We need to compute the coefficients $l_{ij}$ of the braiding.
Since $\rho(x_i)= \overline{K_i}\ot x_i$ It follows that
$$c(x_i\ot x_j)=\overline{K_i}\cdot
x_j\ot x_i=$$
 By (\ref{commute}) $K_i\cdot x_j= \kappa_j^i(\kappa_{j+1}^i)^{-1} x_j$ which is given
 explicitly by:
$$\kappa_j^i(\kappa_{j+1}^i)^{-1}=
\begin{cases}
q_{j,i}q_{j+1,i}\minus \quad &j<i-1,\\
q_{i-1,i}r\minus &j=i-1,\\
rp_{i,i+1}\minus\quad & j=i,\\
p_{i,j}p_{i,j+1}\minus &j>i
\end{cases}
$$
It follows that $\K_i\cdot x_j=K_i\minus K_{i+1}\cdot
x_j=l_{ij}x_j$ where
$$l_{ij}=
\begin{cases}
q_{j,i+1}q_{j+1,i+1}\minus q_{j,i}\minus q_{j+1,i}\quad &j<i-1,\\
q_{i-1,i+1}q_{i,i+1}\minus q_{i-1,i}\minus r      &j=i-1,\\
q_{i,i+1}r\minus r\minus
p_{i,i+1}=r\minus      &j=i,\\
rp_{i+1,i+2}\minus p_{i,i+1}\minus p_{i,i+2}     &j=i+1,\\
p_{i+1,j}p_{i+1,j+1}\minus p_{i,j}\minus p_{i,j+1}      &j>i+1
      \end{cases}$$
      Hence $$l_{ij}l_ {ji}=
      \begin{cases}
      \;\,1\;\;=l_{ii}^0\quad  &|i-j|>2,\\
      r^{-2}=l_{ii}^2 \quad & i=j,\\
      \;\,r\;\;=l_{ii}\minus \quad & |i-j|=1
\end{cases}$$
Thus the braiding is of Cartan type $A_n$ as claimed
\end{pf}
The structure of $B(V)$ depends only on $r.$ It is
finite-dimensional when $r$ is a root of unity by \cite{t} or by
\cite{as}. In fact $B(V)$ is isomorphic to the "positive" part of
either $u_q(sl_n)$ or $U_q(sl_n),$  (where $u_q(\frak{g})$ is
Lusztig's finite-dimensional Hopf algebra derived from
$U_q(\frak{g})$ when $q$ is a root of unity).

The other parameters determine the group $G(U_Q).$  The group will
be infinite if any  of the parameters $p_{i,j},r$ is not a root of
unity. This follows from the definition of the generators $K_i$
given in (\ref{ki}). Thus $B(V)$ may be finite-dimensional while
$SLG$ is infinite. Furthermore,  some of the $K_i$'s may be of
finite order while the others are of infinite order.

It is proved in \cite[Lemma 4.2]{cm} that if $r$ is not a root of
unity then the group $G(B^l)$ ($G(B^r)$) is freely generated by
the $n$ elements $K_i$ ($L_i$).

For the finite-dimensional case we show how  different choices of
the parameters $\{p_{ij}\}$ may provide non-isomorphic groups
$G(B^l)$ and $SLG.$
\begin{theo}\thlabel{family}
Let $Q=\{r\ne 1,\,p_{i,j}\in k\}_{1\le i<j\le n}$ be ${n\choose
2}+1$ non-zero elements in $k,$ let $A(R_Q),\,U_Q$ be the FRT
constructions associated with the Yang-Baxter operator $R_Q$ and
let $V$ and $SLG$ be defined as in (\ref{v}) and (\ref{group}).
Assume $r$ is a root of unity of order $N$ and each $p_{i,j}$ is a
root of unity of order $N_{ij}.$ Then
 \ben
 \item  The groups $G(B^l)$ and thus $SLG$ are finite.
 \item If $\{N, N_{ij}\}_{1\leq
i< j\leq n}$ are relatively prime then each $\K_i$ is of order
$Nm_i$ where $m_i\neq m_l$ for $i\neq l.$ In this case
$B(V)\#kSLG$ is the family defined in \cite[Example 7.27]{as}
 \een\end{theo}
\begin{pf}
1. By  \eqref{deltaoq} and \eqref{ki} we have that
$(K_i)^l(T_j^k)=\gd_{jk}(\kappa_j^i)^l.$ Hence the order of $K_i$
depends on the order of $\{\kappa_j^i,\,j=1,\dots,n\}.$ Now,
$p_{ij}q_{ij}=r$ for all $i<j$ hence if the order of each $p_{ij}$
is $N_{ij}$ then the order of each $q_{ij}$ is ${\rm
lcm}\{N,N_{ij}\}$ and thus the order of each $K_i={\rm
lcm}\{N,N_{ij},\,1\le j\le n\}.$
\medskip

\noin 2. For any $1\le i\le n-1$ set $M_i=\Pi_{j=i+1}^n
N_{ij},\,M_n=1 $ and $m_i=M_{i+1}M_i.$ It follows by part 1 the
order of each $K_i$ is $NM_i$ and hence the order of each $\K_i$
is $Nm_i.$
\end{pf}
\medskip

Set $p_{i,i}=1$ and $p_{j,i}=p_{ij}\minus$ for $i<j.$ Recall
\cite[Th. 3]{ast} that $A(R_Q)$ has a normal group-like element
which is central if and only if $P_l=P_k$ for all $l,k$ where
$$P_l=r^l\prod_{j=1}^n p_{l,j}$$
Since for all $i<j,\,q_{i,j}=r/p_{i,j}$ it is not hard to check
that
 \begin{equation}\label{det}
 P_l=r\prod_{j=1}^{l-1}q_{jl}\prod_{j=l+1}^np_{l,j}=\prod_{j=1}^n\kappa_j^l.\end{equation}
 Set
  \begin{equation}\label{sigma} \gs=K_1\cdots K_n \end{equation}
 We have:
\begin{prop} If the determinant element of $A(R_Q)$ is central then
$\gs$ is a central group-like element of $U_Q.$
\end{prop}
\begin{pf}
By (\ref{ki}) we have that for any $i,\;K_i(T_l^k)$ is nonzero if
and only if $l=k.$ Thus \eqref{deltaoq} implies that for any $u\in
U_Q,\,g\in G(U_Q),$
$$gu(T_l^k)=g(T_l^l)u(T_l^k)\quad\text{and}\quad ug(T_l^k)=u(T_l^k)g(T_k^k).$$
Observe that the element $\gs=K_1\cdots K_n$ satisfies
$\gs(T_l^l)=P_l$ for all $l.$ Thus
 if the determinant is central then $\gs(T_l^l)=\gs(T_k^k)$ for all $l,k.$
 This implies  that for all $u\in U_Q,\,k,l$
$$\gs u(T_l^k)=\gs(T_l^l)u(T_l^k)=u(T_l^k)\gs(T_k^k)=u\gs(T_l^k)$$
Now, for any $a,b\in A(R_Q),\,u\in U_Q$ we have $\gs u(ab)=\sum\gs
u_1(a)\gs u_2(b)$ and $u\gs(ab)=\sum u_1\gs(a)u_2\gs(b).$ So we
can prove by induction on the length of monoms in $A(R_Q)$ that
$u\gs(a)=\gs u(a)$ for all $a\in A(R_Q)$ which proves our claim.
\end{pf}

We consider now some properties of the finite-dimensional Hopf
algebras which are obtained when  the bialgebra $A(R_Q)$ admits a
central determinant.
\begin{theo}\thlabel{slg}
Let $Q=\{r\ne 1,\,p_{i,j}\in k\}_{1\le i<j\le n}$ be roots of
unity, let $A(R_Q),\,U_Q$ be the FRT constructions associated with
the Yang-Baxter operator $R_Q$ and let $V$ and $SLG$ be defined as
in (\ref{v}) and (\ref{group}). Assume that the determinant
element of $A(R)$ is central,  then:
 \ben\item All the group-like elements  $K_i$ have the same order, hence all
the $\K_i$'s have the same order.
 \item If the common order of the $K_i$ is relatively prime to $n$ then
$G(B^l)=\gs\times SLG$ where $\gs$ is the central group-like
defined in (\ref{sigma}). \een
\end{theo}
\begin{pf}
1. By \eqref{det}, the order of each $P_i={\rm
lcm}\{N,N_{ij},\;1\le j\le n\}$ which is also the order of $K_i.$
Thus if the $P_i$'s are all equal then the order of the $K_i$'s
are all equal.
\medskip

\noin 2. Note that $K_iK_j\minus\in SLG$ for all $i\ne j,$ hence
$K_i^{n-1}\prod_{j\ne i}K_j\minus\in SLG$ and so $K_i^n=\gs
K_i^{n-1}\prod_{j\ne i}K_j\minus\in \gs SLG$ for all $i.$ If the
order of each $K_i$ is $m$ and $(m,n)=1$ then this implies that
$K_i\in \gs SLG.$
\end{pf}

As it is known, by letting $r=q^2$ and $p_{ij}=q$ for $i<j$ the
corresponding $U_Q$ is $U_q(gl_n)$ if $q$ is not a root of unity
and $u_q(gl_n)$ otherwise.

The following is an example of a finite-dimensional $U_Q$ so that
the determinant of $A(R_Q)$ is central. But unlike $u_q(sl_n)$ and
$u_q(gl_n)$ the generators $\K_i$ and $K_i$ are not free. This
implies that this Hopf algebra can not be obtained by twisting a
known one.

\begin{exam}\label{7}
Let $n=3$ and let $q$ be a $7^{th}$ root of unity. Let
$$r=q\qquad p_{12}=p_{23}=q_{13}=q^2\qquad q_{23}=q_{12}=p_{13}=q\minus$$
Then the determinant element of $A(R_Q)$ is central by
(\ref{det}). A direct computation using \eqref{ki} shows that
$$K_3=K_1^{-2}K_2^3\quad\text{hence}\quad\K_1^2=\K_2.$$
\end{exam}

In the next example the $n-1$ generators $\K_i$ are free, but not
the $K_i.$
\begin{exam}\label{n}
Given $n,$ let $q$ be a root of unity of order $n+1.$ As for the
$1$-parameter deformation $u_q(gl_n),$ let
$$r=q^2\qquad p_{ij}=q_{ij}=q.$$
A direct computation yields that $\gs=\prod_{i=1}^n K_i=1.$ Thus
$SLG=G(B^l)$ by \thref{family}.4.
\end{exam}

\section{$U_Q $ as a double crossproduct}
In this section we show that for any choice of $\{r\ne
1,p_{ij}\}_{1\le i<j\le n}$ the Hopf algebra $U_Q$ can be
considered as a quotient of a double crossproduct which is
quasitriangular in the finite dimensional case.\medskip

We recall first the definition of the double
crossproduct\cite{r2}. Let $B$ and $H$  be bialgebras so that $B$
is a left $H$-module coalgebra, $H$ is a left $B$-module coalgebra
and certain comparability conditions are satisfied. The double
crossproduct $B\bowtie H$\cite[6.43]{m} is is the tensor product
$B\ot H$ with a coalgebra structure given by the $\gD_B\ot\gD_H.$
The multiplication is defined with respect to the two given module
structures. We omit here the general definition of the product;
instead we will consider the following:

Let $H$ be a Hopf algebra with a bijective antipode and $B$ a sub
Hopf algebra of $(H^0)^{cop}.$ Then $$h\rightharpoonup p=\langle
p_2,h\rangle p_1\quad p\leftharpoonup h=\langle p_1,h\rangle
p_2\quad p\rightharpoonup h=\langle p,h_2\rangle h_1\quad
h\leftharpoonup p=\langle p,h_1\rangle h_2$$ are  defined for all
$h\in H,\,p\in B.$ One can define an $H$-module structure on $B$
and a $B$-module structure on $H$ that satisfy all the necessary
conditions for the double crossproduct. In this case the product
in $B\bowtie H$ is given explicitly by
\begin{equation}\label{double}
(p\bowtie h)(p'\bowtie h')=\sum pp_2\bowtie (S\minus
p'_1\rightharpoonup h\leftharpoonup p'_3)h'
\end{equation}

If $H$ is finite dimensional and $B=(H^*)^{cop}$ the the double
crossproduct is the Drinfeld double $D(H).$\bigskip

We recall also more definitions and results from \cite{t}. The
bialgebra $A(R_Q)$ is endowed with an invertible braiding
$<\,|\,>_R.$ The braiding is given on generators by
$$<T_i^l|T_j^k>=R_{ij}^{kl}$$ where $R_{ij}^{kl}$ are given in
\eqref{rijkl}.

Let $\lambda^+,\,\lambda^-,\,\rho^+,\,\rho^-:A(R_Q)\rightarrow
A(R_Q)^*$ be the following maps:
$$ \gl^+(a)=<a\,|\,->_R\qquad\rho^+(a)=<-\,|\,a>_R$$
$$ \gl^-(a)=<a\,|\,->_{R^{-1}}\quad\rho^-(a)=<-\,|\,a>_{R^{-1}}$$
for all $a\in A(R_Q).$ Recall that $\lambda^+$ is an anti-algebra
and a coagebra map given explicitly by:
 \begin{equation}\label{lt}\lambda^+(T_i^j)=
   \begin{cases}
     (r-1)E_i^j \quad & i>j,\\
     K_i \qquad &i=j,\\
     0 \qquad &\text{otherwise}
  \end{cases}
\end{equation}
and $\rho^+$ is an algebra and an anti-coalgebra map given by
\cite[(6.14)]{t} (after a slight modification) by
 \begin{equation}\label{rt} \rho^+(T_i^j)=
   \begin{cases}
     r^{-2}(r-1)S\minus F_i^j \quad & i<j,\\
     L_i\minus \qquad &i=j,\\
     0 \qquad &\text{otherwise}.
  \end{cases}
\end{equation}

Note that
$$\quad B^l=
 <K_i,K_i\minus,E^i_{i+1}> =\text{the algebra generated by}\,\{{\rm im}\lambda^+,K_i\minus\}$$

Set
$$B^r=<L_i,L_i\minus,F^{i+1}_i> =\text{the algebra generated by}\,\{{\rm im}\rho^+,L_i\}$$

Observe that \eqref{deltae} and \eqref{deltaf} imply that $B^l$
and $B^r$ are Hopf algebras.  Moreover, $U_Q=B^lB^r.$

\begin{lemm}\lelabel{dual}
There exists a Hopf algebra injection  $$\theta :
(B^r)\hookrightarrow (B^l)^{0\,cop}.$$
\end{lemm}
\begin{pf}
Define a map $\hat{\theta} : Im(\rho^+)\rightarrow (Im(\gl^+))^*$
as follows: For $v=\rho^+(b),\;w=\gl^+(a),$ let
\begin{equation}\label{teta} \hat{\theta}(v)(w):= \langle a|b\rangle\end{equation}

>From the definition of $\lambda^+$ and $\rho^+$ it follows that
$\hat{\theta}$ is well  defined, injective, algebra and
anti-coalgebra map. We extend $\hat{\theta}$ step by step:

\noin \underline{Step 1}: For any $1\le j\le n$ we wish to extend
$\hat{\theta}(L_j\minus)$ to an element of $\hom(B^l,k)=(B^l)^*.$
In order to have it we need only to define
$\hat{\theta}(L_j\minus)(K_i\minus)$ for all $1\le i\le n.$ Note
that  \eqref{rt} and the definition of $\hat{\theta}$ imply that
$\hat{\theta}(L_j\minus)(K_i)=\rho^+(T_j^j)(\gl^+(T_i^i))=\langle
T_i^i|T_j^j\rangle =\kappa_j^i.$ Thus define:
$$\hat{\theta}(L_j\minus)(K_i\minus)=(\hat{\theta}(L_j\minus)(K_i))\minus
=(\kappa_j^i)\minus.$$ Observe that  $\hat{\theta}(L_j\minus)$ is
multiplicative as an element of  $(Im(\gl^+))^*,$ therefore the
extension of $\hat{\theta}(L_j\minus)$ defines an element in
$G((B^l)^*).$

\noin\underline{Step 2}: We extend the domain of $\hat{\theta}$
by letting
$$\hat{\theta}(L_j)=(\hat{\theta}(L_j\minus))\minus\in G((B^l)^*).$$
for all $1\le j\le n.$

\noin\underline{Step 3}: In order to extend $\hat{\theta}$ to a
map  $\theta: B^r\rightarrow (B^l)^*$ The only undefined values
are those involving terms of the form
$\theta(F_j^{j+1})(K_i\minus).$ Since By \eqref{deltaf} $S\minus
F_j^{j+1}=-L_{j+1}\minus F_j^{j+1}L_j\minus$ it follows that
$$\theta(F_j^{j+1})(K_i\minus)=\theta(S\minus
F_j^{j+1})(K_i)=-(\theta(L_{j+1}\minus)(K_i))\,
(\theta(F_j^{j+1})(K_i))\,(\theta(L_j)(K_i)),$$  where the right
hand side has been defined in the previous steps. Thus define:
 $$\theta(v)(w)=
\begin{cases}
\hat{\theta}(v)(w) &w\in Im(\gl^+)\\
 \hat{\theta}(S\minus v)(K_i) & w=K_i\minus
\end{cases} $$
and extend $\theta(v)$ to monoms containing $K_i\minus$ in $B^l$
with respect to the coproduct in $(B^r)^{cop}.$ Then $\theta$ is
the desired injection.
 \end{pf}

 Before proving the main theorem of this section we wish to precede with some
 calculations.  By abuse of notations denote $\theta(L_j)$ by $L_j.$ Set

\begin{equation}\label{ei}e_i= (r-1)E_{i+1}^i
 \qquad\qquad f_i= r^{-2}(r-1)F_i^{i+1}\end{equation}

Observe that $e_i=\gl^+(T_{i+1}^i)$ by \eqref{lt} and $S\minus
 f_j=\rho^+(T_i^{i+1})$ by \eqref{rt}.

 By \eqref{rijkl} and \eqref{lt} we have
 \begin{eqnarray}\label{lk}
 L_j\minus(K_i)=&\langle T_i^i|T_j^j\rangle &=\;\kappa_j^i\\
 L_j\minus(e_i)=&\langle T_{i+1}^i|T_j^j\rangle &=\;0\nonumber
 \end{eqnarray}

Since $\theta$ is an anti-coalgebra map it follows that
$S\circ\theta=\theta\circ S\minus,$ and thus
$$
 S\theta(f_j)=\theta(S\minus f_j)=\theta(\rho^+(T_j^{j+1})).$$

Hence \eqref{rijkl} and \eqref{teta} imply that
\begin{eqnarray}\label{sfe}
&&S\theta(f_j)(K_i)=\theta(\rho^+(T_j^{j+1}))(\gl^+(T_{i}^i)=
\langle T_{i}^i|T_j^{j+1}\rangle =R_{i,j}^{j+1,i}=0\\
&&S\theta(f_j)(e_i)=\theta(\rho^+(T_j^{j+1}))(\gl^+(T_{i+1}^i)=
\langle T_{i+1}^i|T_j^{j+1}\rangle
=R_{i+1,j}^{j+1,i}=\gd_{ij}(r-1)\nonumber.\end{eqnarray}

Now, by \eqref{deltae} we have
 $S\minus(e_i)=-K_{i}\minus e_iK_{i+1}\minus$ hence \eqref{commute}
 implies that $S^{-2}(e_i)= K_{i+1}K_{i}\minus e_iK_{i+1}\minus K_i=r\minus
 e_i.$ Therefore,
 \begin{equation}\label{ef}
S\minus \theta(f_j)(e_i)=S\theta(f_j)(S^{-2}e_i)=r\minus
S\theta(f_j)(e_i)=\gd_{ij}r\minus(r-1)
\end{equation}

 It follows that
 \begin{eqnarray}\label{fjsminus}
 \lefteqn{\theta(f_j)(e_i)=S\theta(f_j)(S\minus e_i)=}\\
 &=&-S\theta(f_j)(K_i\minus e_iK_{i+1}\minus)\nonumber\\
 &=&-L_{j}\minus(K_{i}\minus)S\theta(f_j)(e_i)L_{j+1}\minus(K_{i+1}\minus)\nonumber\\
 &&(\text{by applying}
 \;\gD^2(S\theta(f_j))=S\theta^{\ot 3}(\gD^2(f_j)\;\,\text{and by \eqref{lk} or \eqref{sfe}})
 \nonumber\\
 &=&-(\kappa_i  ^{j}\kappa_{i+1}^{j+1})\minus
 \gd_{ij}(r-1)\qquad\quad(\text{by \eqref{lk} and
 \eqref{sfe}})\nonumber\\
 &=&-\gd_{ij}r^{-2}(r-1) \qquad\qquad\qquad(\text{by\eqref{kappa}})\nonumber
 \end{eqnarray}

 We identify $B^r$ with its image  $\theta(B^r)\subset
 (B^l)^{0\,cop}.$ Then $B^r\bowtie B^l$ is defined and we use the above calculations to
 prove:
 \begin{theo}\thlabel{double}
 $U_Q$ is a Hopf-algebra quotient of $B^r\bowtie B^l.$
 \end{theo}

 \begin{pf}
Define $\phi:B^r\bowtie B^l\rightarrow U_Q$ on generators by:
 $$\phi(u\bowtie w)=uw,$$ $u\in B^r,\,w\in B^l.$
Since $U_Q=B^lB^r$ it follows that $\phi$ is surjective. Since the
coproduct in the double is the tensor coproduct  it follows that
$\phi$ is a coalgebra map. Thus we need only to check that $\phi$
is an algebra map.

Since $B^r$ and $B^l$ are contained as Hopf algebras in the double
it follows from the definition of $\phi$ that it is enough to
check multiplicity on generators of the form $(\gep\bowtie
w)(u\bowtie 1).$

We will start from  the relations among the group-like elements:
Observe that \eqref{double} implies that:
$$(\gep\bowtie K_i)(L_j\bowtie 1)=L_j\bowtie (L_j\minus\rightharpoonup
K_i \leftharpoonup L_j)=L_j\bowtie K_i.$$ Thus all group-like
elements in the double commute which is preserved in $U_Q$ and so
$$\phi((\gep\bowtie K_i)(L_j\bowtie 1))=\phi(L_j\bowtie
K_i)=L_jK_i=K_iL_j=\phi(\gep\bowtie K_i)\phi(L_j\bowtie 1).$$

Next we check relations between $K_i$ and $f_j.$ Observe that
\begin{eqnarray*}\lefteqn{(\gep\bowtie K_i)(f_j\bowtie 1)=}\\
&=& f_{j,2}\bowtie (S\minus f_{j,1}\rightharpoonup K_i
\leftharpoonup f_{j,3})\quad\,(\text{by \eqref{double}})\\
&=&L_j\minus(K_i)L_{j+1}(K_i)(f_j\bowtie K_i)\qquad(\text{by
applying}\;\gD^2(f_j)\;\text{and
 by \eqref{sfe}})\\
&=& \kappa_j^i(\kappa_{j+1}^i)\minus  (f_j\bowtie
K_i)\qquad\qquad\quad(\text{by \eqref{lk}})
\end{eqnarray*}
Hence by the definition of $\phi,$
\begin{eqnarray*}
 \lefteqn{\phi((\gep\bowtie K_i)(f_j\bowtie 1))=}\\
&=&(\kappa_{j+1}^i)\minus \kappa_j^i\phi(f_j\bowtie
K_i)\qquad(\text{by
above})\\
&=&(\kappa_{j+1}^i)\minus \kappa_j^i
f_jK_i\\
&=& K_if_j\qquad\qquad\qquad\qquad\quad(\text{by \eqref{ki}})\\
&=&\phi(\gep\bowtie K_i)\phi(f_j\bowtie 1). \end{eqnarray*}

Similarly,
$$(\gep\bowtie e_i)(L_j\bowtie 1)=L_j\bowtie
(L_j\minus\rightharpoonup e_i \leftharpoonup
L_j)=\kappa_i^j(\kappa_{i+1}^j)\minus(L_j\bowtie e_i)$$ which by
\cite[5.21]{t} are the same identities as in $U_Q.$

Consider now the relations between the $\{e_i\}$'s and the
$\{f_j\}$'s. For the sake of convenience let $L_iK_j$ denote the
element $L_i\bowtie K_j= (\gep\bowtie K_j)(L_i\bowtie 1)$ in the
double. We show first that:
 \begin{equation}\label{efd} \kappa^{i+1}_{j}(\gep\bowtie e_i)(f_j\bowtie 1)-\kappa_i^{j+1}
 f_j\bowtie e_i= \gd_{ij}r\minus(r-1)( L_{i}K_{i+1}-L_{i+1}K_i)
 \end{equation}
 Indeed,
 \begin{eqnarray*}
 \lefteqn{(\gep\bowtie e_i)(f_j\bowtie 1)=}\\
 &=& \sum f_{j \,2} \bowtie (S\minus f_{j\,1}\rightharpoonup e_i \leftharpoonup
 f_{j\,3})\qquad(\text{by \eqref{double}})\\
 &=& \sum S\minus f_{j\,1}(e_{i\,3})\; f_{j\,3}(e_{i\,1})\;(f_{j\,2})\bowtie
 e_{i\,2}\\
 &=&  L_{j+1}\minus(K_{i})\, f_j(e_i)\, (L_{j+1}\bowtie
 K_i) + L_{j+1}\minus(K_i)\, L_{j}(K_{i+1})\, (f_j\bowtie e_i)\\
 &+& S\minus(f_j)(e_i)\, L_{j}(K_{i+1})\, (L_{j}\bowtie
 K_{i+1})
 \quad(\text{by applying $\gD^2$ to $e_i$ and $f_j$ and by \eqref{lk} and \eqref{sfe}})\\
 &=& -\gd_{ij}r^{-2}(r-1) \kappa^{i}_{i+1}(L_{i+1}\bowtie K_i) +\kappa^i_{j+1}
 (\kappa^{i+1}_{j})\minus (f_j\bowtie e_i)
 \\
 &+&r\minus(r-1)(\kappa^{i+1}_i)\minus
 (L_{i}\bowtie K_{i+1})\qquad
 (\text{by \eqref{lk},\eqref{ef} and \eqref{fjsminus}})
 \\
 &=& \gd_{ij}\kappa_{i+1}^i r^{-2}(r-1)(L_i\bowtie K_{i+1}-L_{i+1}\bowtie K_{i})+ \kappa_i^{j+1}
 (\kappa^{i+1}_{j})\minus  (f_j\bowtie e_i)\qquad(\text{by
 \eqref{kappa}})
 \end{eqnarray*}

 Since by
\eqref{kappa} $\kappa^{i+1}_i\kappa^i_{i+1}=r$ \eqref{efd}
follows. Now, by \cite[5.23a]{t}the following hold in $U_Q:$
\begin{equation}\label{eft}\kappa^{i+1}_{j}E_{i+1}^iF_j^{j+1} -\kappa^i_{j+1}
F_j^{j+1}E_{i+1}^i=\gd_{ij}r(r-1)\minus(L_{i}K_{i+1}-L_{i+1}K_{i})
\end{equation}

Hence:
\begin{eqnarray*}
\lefteqn{\kappa_{i+1}^{j}\phi(\gep\bowtie e_i)\phi(f_j\bowtie 1))
-\kappa_i^{j+1}  \phi(f_j\bowtie e_i)=}\\
&=& (1-r)^2r^{-2}(\kappa_{i+1}^{j}E_{i+1}^iF_j^{j+1}
-\kappa_i^{j+1}  F_j^{j+1}E_{i+1}^i)\qquad(\text{by \eqref{ei}})\\
 &=&\gd_{ij}r(r-1)\minus (1-r)^2r^{-2}
(L_{i}K_{i+1}-L_{i+1}K_{i})\qquad(\text{by \eqref{eft}})\\
&=&r\minus(r-1)\gd_{ij} (L_{i}K_{i+1}-L_{i+1}K_i)\\
&=&\phi\left(\kappa_{i+1}^{j}(\gep\bowtie e_i)(f_j\bowtie
1)-\kappa_i^{j+1} (f_j\bowtie e_i)\right)\qquad\quad(\text{by
\eqref{efd}})
\end{eqnarray*}
Therefore,
$$ \phi\left((\gep\bowtie e_i)(f_j\bowtie 1)\right)= \phi(\gep\bowtie e_i)\phi(f_j\bowtie
1)).$$ This conclude the proof that $\phi$ is
multiplicative.\end{pf}

 We have shown that $U_Q$ is a homomorphic image of the double
 crossproduct. Note that the relations among the $\{E_i^j,
 F_j^i\}_{i<j}$ are the same in $U_Q$ and in the double. Thus the
 kernel of $\phi$ may contain relations only among the group-like
 elements $\{K_i,L_j\}.$

\begin{coro} If $U_Q$ is finite dimensional then $B^l\cong (B^r)^{*cop}.$
Thus the double is isomorphic to $D(B^l)$ and so $U_Q$ is
quasitriangular.
\end{coro}
\begin{pf}
If $U_Q$ is finite dimensional then  $B^l={\rm im}\gl^+$ and
$B^r={\rm im}\rho^+.$ By \leref{dual} $B^r\subset (B^l)^*$ as
vector spaces. But one can prove similarly that $B^l\subset
(B^r)^*$ via $\phi^*.$ Hence they all have the same dimension and
the double crossproduct is indeed a Drinfeld double.
 \end{pf}\bigskip

 \noin{\bf Acknowledgement} We wish to thank M. Cohen for many
 fruitful discussions and for helpful comments.

\end{document}